\documentclass[12pt]{amsart}
\newtheorem{theorem}{Theorem}[section]
\newtheorem{lemma}[theorem]{Lemma}

\theoremstyle{definition}
\newtheorem{definition}[theorem]{Definition}

\theoremstyle{remark}
\newtheorem{remark}[theorem]{Remark}
\numberwithin{equation}{section}

\def\P{\mathbb P}
\def\E{\mathbb E}
\def\RE{\mathbb R}

\def\l{\lim_{t\uparrow \infty }\,}
\def\p{\par\noindent}

\def\F{{\mathcal F}}

\begin{document}

\title[Flux across hypersurfaces for diffusion processes]
{Asymptotic flux across hypersurfaces for diffusion processes}

\author{Andrea Posilicano}

\address{Dipartimento di Scienze, Universit\`a dell'Insubria, 
via Valleggio 11, I-22100
Como, Italy}

\email{andreap@uninsubria.it}

\author{Stefania Ugolini}
\address{Dipartimento di Matematica, Universit\`a di Matematica, via
Saldini 50, I-20133 Milano, Italy}

\email{ugolini@mat.unimi.it}

\begin{abstract}
We suggest a rigorous definition of the 
pathwise flux across the boundary of a bounded open set for transient 
finite energy diffusion
processes. The expectation of such a flux has the property of
depending only on the current velocity $v$, 
the nonsymmetric (with respect to time reversibility)
part of the drift. In the case where the diffusion has a limiting velocity
we define the asymptotic flux across
subsets of the sphere of radius $R$, when $R$ tends to infinity, and
compute its expectation in
terms of $v$.
\vskip 10pt\noindent
{\it Keywords:} 
diffusion processes, (random) flux across surfaces, current velocity

\end{abstract}
\maketitle

\section{Introduction}

In a previous paper \cite{[PU]} the authors gave a pathwise
probabilistic versions of the
Scattering-into-Cones and Flux--across--Surfaces theorems in Quantum 
Mechanics and then recovered the known analytical results by taking 
suitable expectations.\par
Here we extend the main probabilistic results contained in \cite{[PU]} to
a large class of Markovian diffusions with no a priori connection
with Quantum Mechanics.\par
We consider diffusions on $\RE^d$ which are weak solutions of 
s.d.e.'s of the form
$$dX_t=b(t,X_t)\,dt+dW_t\,, \eqno(1.1)$$ 
where $W_t$ is a standard Wiener process and the drift vector 
field $b_t(x)\equiv b(t,x)$ satisfies
$$
\forall\, T>0,\qquad\E\int_{0}^{T}\|b(t,X_t)\|^2\,dt<+\infty\,. \eqno(1.2)
$$
The condition (1.2) is a {\it finite energy} condition in the sense of
F\"ollmer \cite{[Fo]}. By Girsanov theory one proves that (1.2) is
equivalent to a {\it finite entropy} condition: the probabilitity
measure describing the weak solution of (1.1) has finite relative entropy
with respect to the Wiener measure. By definition of relative entropy,
this fact implies absolute
continuity and therefore the distribution of $X_t$, $t>0$ is
absolutely continuous (with respect to Lebesgue measure) with some density function $\rho_t$.\par
A very peculiar consequence of (1.2) is that not only the Markovian
property is preserved under time reversal, but also the diffusion one.
We will call this property {\it time reversibility of the 
diffusive character}.
The class of diffusion processes individuated by this invariance
property was firstly proposed by Nelson in 1966 \cite{[N1]} in the framework of
Stochastic Mechanics. Such diffusions were then studied by Zheng and Meyer 
who called them {\it semimartingales dans les deux directions du temps} 
\cite{[ZM]}.
Within this class, Carlen in 1984 \cite{[C1]} (see \cite{[Ca]}
for an alternative proof) solved the existence
problem of weak
solutions of (1.1) in the case of (unbounded) drift fields satisfying (1.2).
Successively, F\"ollmer \cite{[Fo]} gave a very elegant characterization  
based on the relative entropy approach. Under time reversal the process
solution of (1.1) is again a solution of a s.d.e. of type (1.1) with some
dual drift field $b^{*}_{t}$ (and, of course, another standard Wiener process) 
which satisfies the relation:
$$b_{t}(x)-b^{*}_{t}(x)=\nabla \log \rho_{t}(x) \eqno(1.3)$$
The duality relation (1.3) allows to introduce a relevant 
decomposition of the drift field as the sum of two vector fields :
$$ v_{t}=\frac{b_{t}+b^{*}_{t}}{2} ,\qquad u_{t}=\frac{b_{t}-b^{*}_{t}}{2} $$
called respectively current and osmotic velocity \cite{[N2]}. 
In the symmetric case, $v_{t}=0$ and thus $b_{t}=u_{t}$ is of 
gradient type according to (1.3). Therefore the current 
velocity $v_{t}$ represents the non symmetric part of the drift
field. We will see that only the current velocity is involved in the 
expression of the flux across surfaces by the diffusion paths.

A useful consequence of the time reversibility of the diffusion 
property is then the validity (in the weak sense) of the continuity
equation for the couple 
$(\rho, v)$:
$$ 
\frac{\partial}{\partial t}\,\rho_t= -\nabla\cdot 
(\rho_t v_t)\,. \eqno(1.4)
$$
In this paper we give a rigorous definition of the pathwise 
flux across the boundary of a bounded set by a transient Markov
diffusion process
solution of (1.1) and satisfying the energy condition (1.2). 

Given a bounded open set $D$ let us consider the function
$$
N_{\partial D}:=
N^+_{\partial D}-N^-_{\partial D}\ ,\eqno(1.5)$$
where $N^+_{\partial D}(\gamma)$ (resp. $N^-_{\partial D}(\gamma) $) 
denotes the number of inward (resp. outward) crossing
by the path $ t\mapsto X_t(\gamma)$ of the boundary ${\partial D}$,
$\gamma$ being the point in the probability space. It is 
not a local time because we need to distinguish outward from inward crossing 
in order to have the pathwise analogoue of a net flux.
The problem is that almost surely the diffusion $X_t$ intersects ${\partial D}$
on a set of times that has no isolated points and is
uncountable. Therefore the 
definition of $N_{\partial D}$ given above makes no sense in
general. However, by a suitable
redefinition of $N_{\partial D}$ as the total mass of the almost surely
compactly supported random distribution 
$-\frac{d}{dt}\,\chi_D(X_t)\ $, where $\chi_D$ is
the characteristic function of the set $D$ (see section 3 for
details), we can give a rigorous definition of the pathwise flux across
$\partial D$.  Then, by using the continuity equation (1.4), in the
case where $\partial D$ is a regular hypersurface we can
compute (see theorem 3.2) the expectation $\Phi_{\partial D}$ of 
the pathwise flux $N_{\partial D}$ in terms 
of the current velocity $v$, obtaining 
$$
\Phi_{\partial D}=\int _{0}^{+\infty}\int_{\partial
D}\rho  (t,x)\,v(t,x)
\cdot n(x)\,
d\sigma(x)\,dt\ ,\qquad\eqno(1.6)
$$
where $n$ denotes the outward unit normal vector along $\partial D$ 
and $\sigma$ is the
surface measure. We
interpret $\Phi_{\partial D}$ as the flux of $X_t$ across $\partial
D$. Note that, by our choice of signs in (1.5), $D$ is a source if 
$\Phi_{\partial D}>0$ and is a sink if $\Phi_{\partial D}<0$.

Of course the definition of flux given above does not extend to the
case of a hypersurface which is not a boundary. This restriction can be
avoided, at least asymptotically, in the case where the diffusion $X_t$ 
has a limiting velocity, i.e.
$$
 \l\frac{1}{t}\, X_t=v_{\infty} \eqno(1.7) 
$$
exists almost surely for some non zero random vector $v_{\infty}$.
A simple but general condition giving (1.7), which is again a 
finite entropy condition, was obtained by Carlen \cite{[C5]} (see Theorem
4.3 below).

Suppose that $\Sigma$ is an open subset of the unit sphere with 
$\partial\Sigma$ a finite union of $C^1$ manifolds. In order to define
the asymptotic flux across $\Sigma$ we consider the cone 
$C_\Sigma:=\left\{\lambda x\in\RE^d\ :\ x\in\Sigma,\ \lambda>
0\right\}$. Using (1.7) (see section 5 for the details) we can then
define $N^a_\Sigma$, the asymptotic pathwise flux
across $\Sigma$, by the limit $R\uparrow \infty$ of
the mass 
$N_{\Sigma_R}$ of the almost surely compactly supported
 random distribution $\frac{d}{dt}\chi_{C\cap
B_R^c}(X_t)$, $B_R$ being the closed ball of radius $R$. Note the change of sign in
the definition of $N_{\Sigma_R}$ with respect to $N_{\partial
D}$. This is consistent with the fact that the exterior normal to
$\Sigma_R:=C_\Sigma\cap S_R$, $S_R$ being the sphere of radius $R$,
coincides with
the interior normal to the boundary of 
$C\cap B_R^c$. No confusion can arise between the two different
definitions since $\Sigma_R$ is never the boundary of an open subset
of $\RE^d$.  

We show (see theorem 5.2) that $N^a_{\Sigma}$ is
well defined since almost surely one has 
$$
N^a_{\Sigma}:=\lim_{R\uparrow\infty}\,N_{\Sigma_R}=\chi_{C_\Sigma}
(v_{\infty})\,.\qquad\eqno(1.8)
$$
Moveover, if as before we define the flux by taking the expectation of
the corresponding pathwise object, we have (see theorem 5.3)
$$
\Phi_\Sigma^a=\lim_{R\uparrow \infty}\,
\int_{0}^{+\infty}\int_{\Sigma_R}\rho(t,x)\,v(t,x)
\cdot n(x)\,
d\sigma(x)\,dt\ .\eqno(1.9)
$$ 
The conditions required in order to obtain the stated
results are, besides (1.1) and (1.2), that
$\rho_t $ and $ v_t$ belong to the Sobolev space $H^1(\RE^3)$ in order to
obtain (1.6), (1.7) to get (1.8), and moreover we assume  
$$
\E\int_{0}^{T}\|\nabla v(t,X_t)\|^2\,dt<+\infty 
$$
for some $T> 0$, to get (1.9).

Finally let us remark that our definition of a pathwise flux across
$\partial D$ is euristically equivalent to the ill-defined Stratonovich
stochastic integral
$$\int _{0}^{+\infty}\left(\int_{\partial D}\delta
(X_t-x)\,n(x)\,d\sigma(x)\right)
\circ 
dX_t \,.
$$
Indeed, proceeding euristically, one obtains 
\begin{align*}
N_{\partial D}=&\int_0^{+\infty}-\frac{d}{dt}\chi_D(X_t)\,dt=
-\int_0^{+\infty}\frac{dX_t}{dt}\cdot\nabla\chi_D(X_t)\,dt\\
=&\int_{0}^{+\infty}\frac{dX_t}{dt}\cdot n(X_t)\,\delta_{\partial D}(X_t)\,dt\\
=&\int_{0}^{+\infty}\left(\int_{\partial D}\delta (X_t-x)\,n(x)\,d\sigma(x)
\right)\circ 
dX_t\,.
\end{align*}

\section{The class of finite energy diffusion processes}

Consider the measurable space $(\Omega,\F)$, 
with $\Omega=C(\RE_+;\RE^d)$, $\F$ the Borel $\sigma$--algebra, and let 
$(\Omega,\F,\F_t,X_t)$ be the evaluation stochastic process 
$X_t(\gamma):=\gamma(t)$, $\gamma\in\Omega$, with $\F_t=\sigma(X_s,\,0\le s\le t)$ the natural 
filtration. 

Let us suppose that:
\smallskip\p
H1) there exists a Borel probability measure $\P$ on 
$(\Omega,\F)$ such that:\p
\smallskip\p
- $(\Omega,\F,\F_t,X_t,\P)$ is a Markov process;
\smallskip\p
- $\displaystyle W_t:=X_t-X_{0}-\int_{0}^tb(s,X_s)\,ds$ is a 
$(\P,\F_t)$--Wiener process, i.e. $\P$ is a weak solution of the 
stochastic differential equation 
$$dX_t=b(t,X_t)\,dt+dW_t \eqno(2.1)$$
with initial 
distribution $\mu = \P\circ X_0^{-1}$, $\mu\ll\lambda$, $\lambda$ the
Lebesgue measure on $\RE^d$;
\smallskip\p
H2) the adapted process $b_t$ satisfies
$$
\forall\,T>0\,,\qquad\E\int_{0}^{T}\|b(t,X_t)\|^2\,dt<+\infty \eqno(2.2)
$$
where $\E$ denotes the expectation with respect to the measure $\P$.

\begin{remark}
Since 
$$
H_{\F_T}(\P,\P^W):=\E\left(\log\left.\frac{d\P}{d\P^W}\right|_{\F_T}\right)
=\frac{1}{2}\,\E\int_{0}^{T}\|b(t,X_t)\|^2\,dt
$$
where $\P^W:=\int \mu (x)\P_{x}^W$, with $\P_{x}^W$ denoting the
Wiener measure starting from $x$, see \cite{[Fo]}, 
the finite energy condition (2.2) is equivalent to a finite relative
entropy condition. Thus (2.2) implies that
$
\P$ is absolutely continuous with respect to
$\P^W$, and so $X_t$ admits a density function $\rho_{t}$ 
for any $t\ge 0$.
\end{remark}
As a consequence of H1 and H2, the Markovian diffusion $X_t$ 
preserves the diffusion property under time reversal. Indeed by
F\"ollmer \cite{[Fo]} one has:
\begin{lemma} Under the hypothesis H1, H2, defining
$\hat\P:=\P\circ R$, where $R$ is the pathwise time reversal on
$C([0,T];\RE^d)$, $R(\gamma)(t):=\gamma(T-t)$, there exists an adapted
process $\hat b_{t}$ such that:
$$\displaystyle \hat W_t:=X_t-X_{0}-\int_{0}^t\hat b(s,X_s)\,ds \eqno(2.3)$$ 
is a 
$(\hat\P,\F_t)$--Wiener process.
\end{lemma}
\begin{remark}
Lemma 2.2 states that the finite energy (entropy) condition (2.2) is a 
sufficient condition for the time reversibility of the diffusion property.
The proof is based on the fact that the finite entropy condition 
is invariant under time reversal. The extension to the infinite
dimensional case is in \cite{[FW]}.
Sufficient conditions are also given in \cite{[HP]}. Sufficient and necessary 
conditions for reversibility of diffusion property, in the case of
Lipschitz drift fields, are investigated in
\cite{[MNS]}.
\end{remark}
It is well known that the drift field can be seen as a stochastic
forward derivative in the sense of Nelson \cite{[N1]}, \cite{[N2]}. In 
particular from (2.1) and (2.2) it follows that (see \cite{[Fo]}):
$$ b_{t}=\lim_{h\downarrow0}\frac{1}{h}
\E[X_{t+h}-X_t|\F_t]\qquad\text{\rm in $L^{2}(\P)$}\,. $$
Analogously, from (2.3) one also has:
$$\hat b_{t}=\lim_{h\downarrow0}\frac{1}{h} \hat \E[X_{t+h}-X_t|\F_t]
\qquad\text{\rm in $L^{2}(\hat\P)$}\,,$$
where $\hat \E$ denotes the expectation with respect to the measure 
$\hat \P$.
For our approach it is convenient to work with the same probability
measure $\P$ as proposed by Nelson \cite{[N1]}. To this end we write:
$$\hat b_{t}=\lim_{h\downarrow0}-\frac{1}{h}
\E[X_{T-t}-X_{T-t-h}|\hat\F_{T-t}]\circ R \qquad\text{\rm in $L^{2}(\P)$}\,,$$
where $\hat\F_t=\sigma(X_s,\, s\ge t)$ is the natural future 
filtration.

Since the Markov property is preserved under time reversal also
the dual drift is given by some measurable function $ \hat
b_{t}(\gamma)=\hat b(t, (X_{t}(\gamma))$. 
Let us define $ b^{*}_{t}(x)=-\hat b_{T-t}(x)$ so
that, as already obtained in \cite{[N2]}:
$$b^{*}_{t}=\lim_{h\downarrow0}-\frac{1}{h}\E[X_{t}-X_{t-h}|\hat\F_{t}] 
\qquad\text{\rm in $L^{2}(\P)$}\,,$$
and the following relation holds:
$$b_{t}(x)-b^{*}_{t}(x)=\nabla \log \rho_{t}(x) \eqno(2.4)$$
between the drift field and its dual (see \cite{[C1]} and for the non 
Markovian case \cite{[Fo]}).

The duality relation (2.4) allows to introduce the  
decomposition 
$$ b_{t}=u_{t}+v_{t},\qquad b^{*}_{t}=-u_{t}+v_{t}$$
where 
$$ v_{t}=\frac{b_{t}+b^{*}_{t}}{2} ,\qquad u_{t}=\frac{b_{t}-b^{*}_{t}}{2} $$
are called current and osmotic velocity respectively \cite{[N1]}.
 
In the symmetric case, $v_{t}=0, b_{t}=u_{t}, b^{*}_{t}=-u_{t}$, 
thus the drift field coincides up to the sign with its dual and, 
according to (2.4), it is of 
gradient type. Therefore the current 
velocity $v_{t}$ represents the not symmetric (with respect to time
reversal) part of the drift
field. We will see that only the current velocity is involved in the 
expression of the flux across surfaces.

An important consequence of the time reversibility of the diffusion 
property is the validity of the continuity equation (in the weak
sense) for the couple 
$(\rho, v)$.

Indeed, recalling the Fokker-Planck equation associated with (2.1):
$$ 
\frac{\partial}{\partial t}\,\rho_t= -\nabla\cdot 
(\rho_t b_t)+ \Delta\rho_t
$$
and the Fokker-Planck equation associated with (2.3):
$$ 
\frac{\partial}{\partial t}\,\rho_t= -\nabla\cdot 
(\rho_t b_t^{*})- \Delta\rho_t
$$
and putting together the two equations one has:
$$ 
\frac{\partial}{\partial t}\,\rho_t= -\nabla\cdot 
(\rho_t  v_t) \eqno(2.5)
$$
where the definition of the current velocity has been used.

\section{The pathwise flux across a boundary.}

Given an open set $D$ we want now to define the 
flux across $\partial D$ by the path of a diffusion.

In order to do this we would like to introduce a pathwise analogous of
the flux
as the function 
$$
N_{\partial D}(\gamma):=
N^+_{\partial D}(\gamma)-N^-_{\partial D}(\gamma)\ ,
$$
where $N^+_{\partial D}(\gamma)$ (resp. $N^-_{\partial D}(\gamma) $) 
denotes the number of inward (resp. outward) crossing
by $[0,+\infty)\ni t\mapsto\gamma(t)$ of $\partial D$.
The
problem is that the above definition makes no sense since
$\P$--a.s. the set $\{t\, :\, X_t\in \partial D\}$ has no isolated points and 
is uncountable. Therefore we are forced to proceed in an alternative
way.\par
Let us observe that if 
$\#\left\{t\,:\,\gamma(t)\in \partial D\right\}<+\infty$ then
$N_{\partial D}(\gamma)$ is the total mass of
the random signed measure 
$$
\sum_{t\in\{s\,:\,\gamma(s)\in \partial D\}}c(t)\,\delta_{t}\ ,
$$
where $c(t)=+1$ if $t$ corresponds to an
outward crossing and $c(t)=-1$ if $t$ corresponds to an
inward crossing. Therefore
$$
\sum_{t\in\{s\,:\,\gamma(s)\in \partial D\}}c(t)\,\delta_{t}=-
\frac{d}{dt}\, \chi_{D}(\gamma(t))\,,
$$
where the derivative has to be intended in distributional sense, and
thus we give the following
\begin{definition} 
Given any open domain $D$, we define the random distribution 
$$\mu_D:\Omega\to {\mathcal D}'(\RE)$$ 
by
$$
\langle\mu_D(\gamma),\phi\rangle
:=\chi_{D}(\gamma(0))\,\phi(0)+\int_{0}^{+\infty}
\chi_D(\gamma(t))\,\dot\phi(t)\, dt\,,
$$
where $\phi\in {\mathcal D}(\RE)\equiv C_c^\infty(\RE)$.
\end{definition}
Supposing now that $D$ is bounded and that $X_t$ is transient, we have
that $\P$-almost surely the random distribution 
$\mu_{D}$ has a compact support and so
its mass, which we denote by $N_{\partial D}$, is well defined. 
We define then the flux across $\partial D$ by $\Phi_{\partial D}:=
\E(N_{\partial D})$. 

By the continuity equation (2.5) such an expectation
can be explicitly calculated in terms of the current velocity $v$ (use
\cite{[PU]}, theorem 7):
\begin{theorem} Let $(\Omega,\F,\P_t,X_t,\P)$ satify H1 and H2, with 
$\rho_t \in
H^1(\RE^3)$ and $v_t \in
H^1(\RE^3)$ for any $t \ge 0$. 
 For any open bounded
domain $D$, with $\partial D$ a finite union of $C^1$ manifolds, one has
$$
\Phi_{\partial D}=\int _{0}^{+\infty}\int_{\partial
D}\rho  (t,x)\,v(t,x)
\cdot n(x)\,
d\sigma(x)\,dt\ ,
$$
where $n$ denotes the outward unit normal vector along $\partial D$ 
and $\sigma$ is the
surface measure.
\end{theorem}

\section{Diffusion with an asymptotic velocity.}

Since our goal is to define an asymptotic flux across hypersurfaces, we 
need to impose a condition on the time evolution of the process $
\frac{1}{t}\,X_t$.
\begin{definition}
We say that the diffusion paths admit an asymptotic velocity when\p
H3)
$$ 
\l\frac{1}{t}\, X_t=v_{\infty}\ne 0\qquad 
\P\text{--a.s.}
$$ 
and moreover $\mu_{\infty}\ll\lambda$, where
$\lambda$ denotes the Lebesgue measure on $\RE^d$ and $\mu_\infty$ is the
distribution of $v_\infty$.
\end{definition}
From now on by an open cone $C_\Sigma$ we will mean a
set of the form $$\left\{\lambda x\in\RE^3\ :\ x\in\Sigma,\ \lambda>
0\right\}\,,$$ where $\Sigma$ is an open subset of the unit sphere with 
$\partial\Sigma$ a finite union of $C^1$ manifolds.
\begin{remark} For any open cone $C_\Sigma$, any ball $B_R$ of radius 
$R$, and for any diffusion $(\Omega,\F,\P_t,X_t,\P)$ satifying H1, H2
and H3, one has $$
\l \chi_{C_\Sigma\cap B_R^c}( X_t)=\l \chi_{C_\Sigma}( X_t)
=\chi_{C_\Sigma}(v_{\infty} )
\qquad 
\P\text{--a.s.}\ .
$$
See \cite{[PU]} for a two-line proof. Thus hypotheses H3 
requires that the limiting velocity is non negligible 
and such 
that asymptotically the paths have the same direction as their
limiting velocity. For example 
the Brownian motion in $\RE^{3}$ is transient but it has no limiting
velocity according to our definition because of the S.L.L.N. Only
a Brownian motion with drift could satisfy the requirement of our
definition.
\end{remark}
A simple but general condition giving the existence of a limiting
velocity, 
which is again a 
finite relative entropy condition (now on the full $\sigma$-algebra
and with respect to $\widetilde\P$, see the proof below), is given in
the 
following (see \cite{[C5]}):

\begin{theorem} Let $(\Omega,\F,\P_t,X_t,\P)$ satify H1 and H2. 
If moreover one has: 
$$
\E\int_{t_0}^{+\infty}\|b(t,X_t)-X_t/t\|^2\,dt<+\infty\,\qquad t_0>0, 
\eqno(4.1)
$$
then:
$$
 \l\frac{1}{t}\, X_t=v_{\infty}\qquad \P\text{\rm --a.s.}\
$$
for some random variable $v_{\infty}$.
\end{theorem}
\begin{proof}
The condition (4.1) implies, by \cite{[E]}, prop. 2.11, that $\P\ll\widetilde
\P$ on $\sigma(X_s, t_0\le s< +\infty)$, where 
 $\widetilde \P$ is the weak solution of the 
simple stochastic differential equation
$$
dX_t=\frac{1}{t}\,X_t\ dt +\ d\widetilde W_t\ .
$$
with $\widetilde W_t$ a standard Wiener process and is such that
$\P\circ X_{t_0}^{-1}=\widetilde\P\circ X_{t_0}^{-1}$.
Therefore:
$$
d\left(\frac{1}{t}\,X_t\right)=-\frac{1}{t^{2}}\,X_t\ dt+\frac{1}{t}\,dX_t =
\frac{1}{t}\ \widetilde W_t\ .
$$
and so
$$
\frac{1}{t}\,X_t=\frac{1}{t_{0}}\,X_{t_{0}}+\int_{t_0}^{t}\frac{1}{s}\ d\widetilde W_s\ .
$$
Since
$$
\widetilde\E\left(\int_{t_0}^{+\infty}\frac{1}{s}\ d\widetilde 
W_s\right)^{2}=
\int_{t_0}^{+\infty}\frac{1}{s^2}\ ds<+\infty 
$$
by Doob's martingale convergence theorem one gets $\widetilde \P\text{\rm 
--a.s.}$ convergence of $\frac{1}{t}\,X_t$. The theorem then follows by
absolutely continuity.
\end{proof}
\begin{remark} Under the same hypotheses of Theorem 4.3
it is possible to prove (see \cite{[C5]}) that the random variable 
$v_{\infty}$
generates the tail $\sigma$--algebra 
$$
{\mathcal T}:=\bigcap_{t>0}\sigma(X_s,\ s\ge t)\ .
$$ 
\end{remark}

\section{The asymptotic flux across hypersurfaces.}

Let us consider hypersurfaces of the following type:
$$\Sigma_R =C_\Sigma\cap S_R\,,$$
where $S_R$ is the sphere of radius $R$. We will define the pathwise 
flux across $\Sigma_R$ in the limit when $R\uparrow\infty$.\par
Suppose at first that $\#\left\{t\,:\,\gamma(t)\in \Sigma_R\right\}<+\infty$. 
Since, by H3, $t\mapsto\gamma(t)$ is definitively either in
$C_\Sigma$ or in $\bar C_\Sigma^c$, if $R$ is sufficiently large one has 
$$
\sum_{t\in\left\{s\,:\,\gamma(s)\in \Sigma_R\right\}}c(t)\,\delta_{t}=
\sum_{t\in\left\{s\,:\,\gamma(s)\in \Sigma_R
\cup(\partial C_\Sigma\cup B_R^c)\right\}}c(t)\,\delta_{t}
=\frac{d}{dt}\, \chi_{C_\Sigma\cap B_R^c}(\gamma(t))\,.
$$
We are therefore lead to give the following
\begin{definition} The asymptotic pathwise flux across
$\Sigma$ is defined by
$$
N^a_\Sigma:=\lim_{R\uparrow\infty} N_{\Sigma_R}\,,
$$
where $N_{\Sigma_R}$ is the total mass of the random distribution 
$-\mu_{C_\Sigma \cap  B_R^c}$.
\end{definition}
The following result shows that the above definition makes sense.
 
\begin{theorem} Let $(\Omega,\F,\P_t,X_t,\P)$ satify H1, H2 and H3. 
Then $\P$-almost surely the random distribution 
$\mu_{C_\Sigma\cap B_R^c}$ has a compact support and so
its mass $N_{\Sigma_R}$ is well defined. Moreover
one has $$
\lim_{R\uparrow\infty} N_{\Sigma_R}
=\chi_{C_\Sigma}(v_{\infty})\qquad 
\P\text{--a.s.}\, .
$$
\end{theorem}
\begin{proof} Let
$$
\tau_R(\gamma):=\sup\left\{t\ge 0\ :\
\|\gamma(t)\|<R\right\}\ .
$$
By H3, $\tau_R<+\infty$ $\P$--a.s.  
Thus $\mu_{C\cap 
B_R^c}$ has compact support $\P\text{--a.s.}$.\par
Let $\phi_\gamma\in {\mathcal D}(\RE)$ such that
$\phi_\gamma=1$ on a neighbourhood of $[0,\tau_R(\gamma)]$. 
By the definition of $\mu_{C_\Sigma\cap 
B_R^c}$ one has 
\begin{align*}
\langle\mu_{C_\Sigma\cap 
B_R^c}(\gamma),\phi_\gamma\rangle
=&-\chi_{{C_\Sigma\cap 
B_R^c}}(\gamma(0))-\chi_{C_\Sigma}(v_{\infty}(\gamma))\int_{\tau_R(\gamma)}^{+\infty}\dot\phi_\gamma(t)\,
dt\\
=&-\chi_{{C_\Sigma\cap 
B_R^c}}(\gamma(0))+\chi_{C_\Sigma}(v_{\infty}(\gamma))\ ,
\end{align*}
and the thesis then immediately follows by taking the limit $R\uparrow\infty$.
\end{proof}
The next theorem shows that the definition of asymptotic flux across
$\Sigma$ by $\Phi^a_\Sigma:=\E(N^a_\Sigma)$ is consistent with the 
result given in theorem 3.2 in the case of the flux across a boundary:  
\begin{theorem} Let $(\Omega,\F,\P_t,X_t,\P)$ satify H1, H2 and H3 and
suppose $\rho_t\in
H^1(\RE^3)$, $v_t \in
H^1(\RE^3)$ for any $t \ge 0$ and 
$$
\E\int_{0}^{T}\|\nabla v(t,X_t)\|^2\,dt<+\infty \eqno(5.1)
$$
for some $T> 0$. Then 
$$\Phi^a_\Sigma\equiv\E(\chi_{C_\Sigma}(v_\infty))=
\lim_{R\uparrow \infty}\,
\int_{0}^{+\infty}\int_{\Sigma_R}\rho(t,x)\,v(t,x)
\cdot n(x)\,
d\sigma(x)\,dt\ .
$$
\end{theorem}
\begin{proof} Proceeding as is \cite{[PU]} one has 
$$
\E(\chi_{C_\Sigma}(v_{\infty}))=
\lim_{R\uparrow \infty}\,
\int_{0}^{+\infty}\int_{\Sigma_R\cup ( \partial C_\Sigma\cap B_R^c)}\rho(t,x)\,v(t,x)
\cdot n(x)\,
d\sigma(x)\,dt\ .
$$
The proof is then concluded by proving that: 
$$
\lim_{R\uparrow \infty}\,
\int_{0}^{+\infty}\int_{ \partial C_\Sigma\cap B_R^c}
\rho(t,x)\,v(t,x)
\cdot n(x)\,
d\sigma(x)\,dt\ =0\, .\eqno(5.2)
$$
By the monotone convergence theorem, (5.2) follows from
$$
\int_{0}^{T}\int_{ \partial(C_\Sigma\cap B_R^c)}
|\rho(t,x)\,v(t,x)
\cdot n(x)|\,
d\sigma(x)\ dt\ <+\infty \eqno(5.3)
$$
for some $T> 0$.
Since
\begin{align*}
&\int_{0}^{T}\int_{ \partial(C_\Sigma\cap B_R^c)}
|\rho(t,x)^{1/2}\rho(t,x)^{1/2}\,v(t,x)
\cdot n(x)|\,
d\sigma(x)\,dt\\ 
\le& \int_{0}^{T}\left(\left(\int_{ \partial(C_\Sigma\cap B_R^c)}
\rho(t,x)\,d\sigma(x)\right)^{1/2}\right.\times\\
&\left.\left(\int_{ \partial(C_\Sigma\cap
B_R^c)}\rho(t,x))\,\|v(t,x)\|^2\,d\sigma(x)\right)^{1/2}\right)\,dt
\,, 
\end{align*}
by trace estimates on functions in $H^1(\RE^3)$ (see e.g.\cite{[B]}, chap.5)
one has
\begin{align*}
&\int_{0}^{T}\int_{ \partial(C_\Sigma\cap B_R^c)}
|\rho(t,x)\,v(t,x)
\cdot n(x)|\,
d\sigma(x)\,dt\\ &
\le \int_{0}^{T}\left(\int_{ \RE^3}\left(
\rho(t,x)\,dx +  \int_{ \RE^3}\|\nabla \rho^{1/2}(t,x)\|^2 dx
\right)^{1/2}\right.
\times\\
&\left.\times\left(\int_{ \RE^3}
\rho(t,x) \|v(t,x)\|^2\,dx +  \int_{ \RE^3}\|\nabla (\rho^{1/2} v)(t,x)\|^2\,dx \right)^{1/2}\right)\,dt\,.   
\end{align*}
From (2.4) we have $ u=\frac{{\nabla \rho }}{{2 \rho}}$, hence one has
$\nabla \rho^{1/2} = u \rho^{1/2}$ and $\nabla (\rho^{1/2}v) = u v \rho^{1/2}+
\rho^{1/2} \nabla v$. Therefore in order to obtain (5.3) it is  
sufficient to have (5.1) and
$$
\E\int_{0}^{T}\left(\|u(t,X_t)\|^2\ + \|v(t,X_t)\|^2\right)\,dt<+\infty 
$$
which is equivalent 
to (2.2).  
\end{proof}
\section{conclusion} The present paper introduces the notion of 
flux across the boundary of a bounded open set $D$ in Euclidean
space. For regular compact domains it is the expectation of the total
mass of the (almost surely compactly supported) random distribution
given by the distributional time-derivative of the functional 
$$
-\chi_D(X_t)\,,
$$
where the diffusion $X_t$, which satisfies the stochastic differential
equation $dX_t=b_t(X_t)\,dt+dW_t$, is transient. If, moreover, the limit 
$$
\lim_{t\uparrow\infty}\, \frac{1}{t}\,X_t
$$
exists almost surely, then to the notion of asymptotic flux is given a
sense as well. The flux can be expressed in terms of $\rho_t$, the
density of the distribution 
of $X_t$, and the current velocity
$v_t=b_t-\frac{1}{2}\,\nabla \log\rho_t$.

\end{document}